\documentclass[a4paper,12pt]{article}
\usepackage{amsmath,amsfonts,amssymb}

\def\l{\lambda}

\def\x{\mathbf x}
\def\y{\mathbf y}
\def\z{\mathbf z}
\def\N{{\mathbb N}}
\def\Z{{\mathbb Z}}

\def\mfS{{\mathfrak S}}

\def\moins{\raise 1pt\hbox{{$\scriptstyle -$}}}
\def\plus{\!\raise 1pt\hbox{{$\scriptstyle +$}} }

\newtheorem{theorem}{Theorem}
\newtheorem{proposition}[theorem]{Proposition}

\begin{document}

\title{\bf Generalisation of Scott permanent identity}

\author{Alain Lascoux\thanks{partially supported by the ANR project MARS (BLAN06-2 134516).}\\
\small IGM, Universit\'e de Paris-Est\\
\small 77454 Marne-la-Vall\'ee CEDEX 2\\
\small Alain.Lascoux@univ-mlv.fr\\
 \small {\tt http://phalanstere.univ-mlv.fr/$\sim$al}
} 

\date{}

\maketitle

\begin{abstract}
Let $\x=\{ x_1,\ldots, x_r\}$,
 $\y=\{y_1,\ldots, y_n  \}$,
 $\z=\{z_1,\ldots, z_n  \}$ be three sets of indeterminates.
We give the value of the determinant
$$ \Bigg| \prod_{x\in \x} (xy-z)^{-1} \Bigg|_{y\in\y, z\in\z}   $$
when specializing $\y$ and $\z$ to the set of roots of
$y^n-1$ and $z^n- \xi^n$ respectively.
\end{abstract}

In the case where $r=2$ and  $\x=\{1,1\}$  the determinant
$\big| (y-z)^{-2}   \big|_{y\in\y, z\in\z}$    factorizes 
into the determinant  of the Cauchy matrix $[ (y-z)^{-1}]$ and its 
permanent. Scott \cite{Scott,Minc}  found the value of this permanent
when specializing $\y$ to the roots of $y^n-1$ and $\z$ to the roots
of $z^n+1$. Han \cite{GuoNiu} described more generally the case where
$\z$ is the set of roots of $z^n + az^k+b$ instead of $z^n+1$.

Instead of restricting to $r=2$ and specializing $\x$, 
we shall consider the determinant 
$$ \Bigg| \prod_{x\in \x} (xy-z)^{-1} \Bigg|_{y\in\y, z\in\z}  \, , $$
and obtain in Th. \ref{ThGaudinScott} 
its value when specializing $\y$ and $\z$. 
The remarkable feature is that this value is a product
of sums of monomial functions in $\x$ without multiplicities,
thus extending the factorized expressions of \cite{Scott,GuoNiu}.

We first need a few generalities about symmetric functions \cite{Cbms}.

Given two sets of indeterminates $\x,\z$ (we say alphabets),
the complete functions $S_n(\x-\z)$ are the coefficients of
the generating function
$$ \sum_n \gamma^n S_n(\x-\z)= \prod_{z\in\z}( 1-\gamma z)
    \prod_{x\in\x}( 1-\gamma x)^{-1} \, .$$

For any $r$, any $\l\in\Z^r$, $S_\l(\x-\z) = 
  \det\bigl( S_{\l_i+j-i}(\x-\z)\bigr)$.

In the case where $\z=0$, and $\x$ of cardinality $r$, these functions
can be obtained by symmetrisation over the symmetric group
$\mfS_r$.       Let $\pi_\omega$ be the following 
operator on functions  in $\x$:
$$  f\to f\pi_\omega := \sum_{\sigma\in\mfS_r} 
  \left( f\, \prod_{1\leq i<j\leq r} (1-x_j/x_i)^{-1} \right)^\sigma \, . $$
Then, when $\l\geq [1\moins r,\ldots, \moins 1,0]$ (i.e. $\l_1\geq 1\moins r$,
 \ldots, $\l_r\geq 0$),
the monomial $x^\l=x_1^{\l_1}\cdots x_r^{\l_r}$ is sent onto 
$S_\l(\x)$ under $\pi_\omega$. When $\l$ is a partition (i.e. 
$\l_1\leq \cdots \l_r\leq 0$), $S_\l(\x)$ is the Schur function of index $\l$.

\bigskip
Let $n$ be a positve integer, 
$\xi$ and indeterminate and $\z$ be the set of roots
of $z^n-\xi^n$. Equivalently, $e_i(\z)=0$ for $1 \leq i \leq n\moins1$,
 $e_n(\z)= (\moins1)^{n-1} \xi^n$. 
For any integer $j$, any $\x$, one has
$$ S_j(\x-\z)= S_j(\x) - \xi^n S_{j-n} (\x)   \, ,$$
and more generally, from the determinantal expression of Schur functions,
for any $\l\in\N^r$,
$$ S_\l(\x-\z)= \sum_{u\in \{0,n\}^r} S_{\l-u} (\x)  (-1)^{|u|/n} \xi^{|u|} \, .$$

In particular, when $\l=\underbrace{n\moins 1,\ldots, n\moins 1}_{r-1},p$
then the terms with negative last component $n-p$ vanish, and  
 the set  $\big\{ \l-u \big\}$ to consider is $ \big\{[v,p] :\, 
                v\in \{n\moins 1, \moins 1\}^{r-1}\big\}$.
Reordering indices, putting $q=p\moins r\plus 1$, one rewrites the sum as
\begin{equation}   \label{SchurSpec}
 S_\l(\x-\z)=  (-1)^{r-1} \sum_{u\in \{0,n\}^{r-1}} S_{q,u} (\x)
   (-1)^{(|\l| -|u| -q)/n}     \xi^{|\l| -|u| -q}  \, . 
\end{equation}

Since $x^v\pi_\omega= S_v(\x)$, for any 
$v\geq [1\moins r,\ldots, \moins 1,0]$,
one can rewrite (\ref{SchurSpec}) as a symmetrisation of monomial functions
in $\x-x_1=\{x_2,\ldots, x_r\}$~:
\begin{equation}
 S_\l(\x-\z)= (-1)^{r-1}\sum_{j=0}^{r-1} x_1^q m_{n^j}(\x-x_1) 
      (-1)^{r-1-j}   \xi^{|\l|-jn-q} \, \pi_\omega \, .
\end{equation}

From the identity 
$$ m_{n^j}(\x-x_1) = m_{n^j}(\x) -x_1^n m_{n^{j-1}}(\x) 
   + x_1^{2n} m_{n^{j-2}}(\x)  +\cdots + (-x_1^n)^j \, ,$$
one sees that $x_1^q m_{n^j}(\x-x_1)\, \pi_\omega$ is equal to

$$S_q(\x)  m_{n^j}(\x) -S_{q+n} (\x) m_{n^{j-1}}(\x)   +\cdots  
                                            +(-1)^j S_{q+jn}(\x) \, .$$
Since  on the other hand $S_{(n-1)^{r-1}p}(\x-\z)$ belongs to the linear
span of Schur functions, or monomial functions, indexed by
partitions $\mu$ such that $\mu_1\leq n\moins 1$, one can restrict
this last  sum to the term $(-1)^j S_{q+jn}(\x)$.

In summary, one has the following expression for the specialisation
of the Schur function that we are considering.

\begin{proposition}
Let $\x$ be an alphabet of cardinality $r$, $\z$ be the set of roots
$z^n- \xi^n=0$, $p\leq n\moins 1$, $N=(n\moins 1)(r\moins 1)$.  Then
\begin{equation}   \label{SchurSpec2}
 S_{(n-1)^{r-1}p}(\x-\z) = \sum_\mu m_\mu(\x)\, \xi^{N+p-|\mu|} \, ,
\end{equation}
sum over all partitions $\mu\in\N^r$, $\mu_1\leq n\moins 1$.

\end{proposition}

For example, for $n=4$, $r=2$, one has
$$ S_{30}(\x-\z)= m_3(\x)+m_{21}(\x),\ S_{31}(\x-\z)= 
m_{31}(\x)+m_{22}(\x) +\xi^4 \, ,$$
$$ S_{32}(\x-\z)= m_{32}(\x) +\xi^4 m_1(\x), \ S_{33}(\x-\z)=
   m_{33}(\x) +\xi^4(m_{2}(\x)+m_{11}(\x))  \, .$$

\bigskip
Let 
$$ D(\x,\y,\z)=  \Big| \prod_{x\in\x} (xy-z)^{-1} \Big|_{y\in\y,z\in\z} \, .$$

In the case $r=2$, this determinant has been obtained by Izergin and Korepin
\cite{Izergin} 
as the partition function of the Heisenberg XXZ-antiferromagnetic model.  
Gaudin \cite{Gaudin} had previously described the partition function
of some other model as the determinant 
$\big|(x-y)^{-1} (x-y+\gamma)^{-1}   \big|$ for some parameter $\gamma$.  

The Izergin-Korepin determinant  is used in the enumeration of 
alternating sign matrices \cite{Bressoud}. In that case, one first 
specializes $\x= \{ e^{2i\pi/3}, e^{4i\pi/3} \}$. Okada \cite{Okada} 
evaluates more general partition functions corresponding to 
similar determinants or Pfaffians, and to other roots of unity
(see also \cite[Th. 7.2]{Pfaff}).

We shall take another point of view, keep $\x$ generic, 
but specialize instead $\y$ and $\z$.
In \cite[Formula 4]{BigGaudin}, it is shown that the function
$$ G(\x,\y,\z)= \frac{ D(\x,\y,\z)}{\Delta(\z)} \, 
   \prod_{x\in \x} \prod_{y\in \y}\prod_{z\in \z}
    (xy-z)  $$
is equal to the determinant of the matrix
\begin{equation}   \label{MatGaudin}
 \Bigl[  S_{\Box\, j}(y_i\x -\z) 
           \Bigr]_{j=0\ldots n-1,\, i=1\dots n}  \, , 
\end{equation}
where $\Box=(n-1)^{r-1}$, and $\Delta(\z)=\prod_{1\leq i<j\leq n} (z_i-z_j)$.

For any $k\in\N$, let $\varphi_k$ be the sum of all monomial
functions $m_\mu(\x)$ of degree $k$, with $\mu_1\leq n\moins 1$ 
(notice that $\varphi_k=0$ when $k>(n-1)r$).
{}From (\ref{SchurSpec2}), one has that $S_{\Box\, j}(y_i\x -\z)$
specializes, when $\z$ is the set of roots of $z^n-\xi^n$, into
\begin{equation}   \label{SchurSpec3}
 S_{\Box\, j}(y_i\x -\z)= y_i^{N+j} \varphi_{N+j} +
  \xi^n y_i^{N+j-n} \varphi_{N+j-n} + 
  \xi^{2n} y_i^{N+j-2n} \varphi_{N+j-2n} +\cdots   \, .
\end{equation}

Specializing further $\y$ into the roots of $y^n-1$, one sees that 
the matrix (\ref{MatGaudin}) factorizes into the product of the matrix
$\Bigl[ y_i^{(N+j)}   \Bigr]$, where $(k)= k \ \mod n$,
and the diagonal matrix
\begin{multline*}
 diag\Big(
 (\varphi_N +\xi^n \varphi_{N-n} +\xi^{2n} \varphi_{N-2n} +\cdots),
 (\varphi_{N+1} +\xi^n \varphi_{N+1-n} +\xi^{2n} \varphi_{N+1-2n} +\cdots), \\
  \cdots, \  
 (\varphi_{N+n-1} +\xi^n \varphi_{N+n-1-n} +\xi^{2n} \varphi_{N+1-n} +\cdots)
 \Big) \, .
\end{multline*}

For example, for $n=3$, $r=3$, 
\begin{multline*}
 S_{220}(y_i\x-\z)= y_i^4\varphi_4 +y_i\varphi_1 \xi^3, 
 S_{221}(y_i\x-\z)= y_i^5\varphi_5 +y_i^2\varphi_2 \xi^3,   \\ 
 S_{222}(y_i\x-\z)= y_i^6\varphi_6 +y_i^3\varphi_3\xi^3 +\xi^6 \, ,
\end{multline*}
and the matrix factorizes into
$$ \begin{bmatrix}  y_1 & y_1^2 & 1\\ y_2 & y_2^2 & 1\\ y_3 & y_3^2 & 1\\
  \end{bmatrix}  \ 
 \begin{bmatrix} \varphi_4+\varphi_1\xi^3 & 0 &0 \\ 
        0 & \varphi_5 +\varphi_2\xi^3 & 0 \\
       0 & 0 & \varphi_6 +\varphi_3\xi^3 +\xi^6 
 \end{bmatrix}
$$

Taking into account that $\prod_{x,y,z} (xy-z)$ specializes into
 $\prod_{y,x}(x^ny^n -\xi^n) = \prod_{x\in\x} (x^n-\xi^n)^n$,
and that the determinant of powers of the $y\in\y$ is a permutation of
the Vandermonde in $\y$, one obtains the following theorem.

\begin{theorem}    \label{ThGaudinScott}
Let $n,r$ be two positive integers, $N=(n-1)(r-1)$.
 Let $\x$ be an alphabet of cardinality $r$, 
$\y$ be the set of roots of $y^n-1$,
 $\z$ be the set of roots of $z^n-\xi^n$.  Then 
$$ \Delta(\y) \Delta(\z)  \left|  \prod_{k=1}^r (x_k y_i -z_j)^{-1}
    \right|_{i,j=1\ldots n} =  
 \frac{ (-1)^{(n-1)(n/2 +r-1)}}{\prod_{x\in\x} (x^n-\xi^n)^n} \
 \prod_{i=0}^{n-1} \left(\sum_{j=0}^\infty \varphi_{N+i-nj}\,  \xi^{nj}
       \right) \, . $$
\end{theorem}

For  $\x=\{1,1\}$, this theorem is due to Han\cite{GuoNiu}.
In that case, $\varphi_i=i+1$ and $\varphi_{n-1+i}=n-i$
for $i=0,\ldots, n\moins 1$, and the product appearing in the theorem is
$$ n\, (n-1+\xi^n)\, (n-2+2\xi^n)\cdots (1+(n-1)\xi^n)  \, .    $$

For $r=5$, $n=3$, as a further example, the theorem furnishes the expression  
\begin{multline*}
 \prod_{k=1}^5  (x_k^3-\xi^3)^{-3} \,
 \left( \varphi_8 +\varphi_5\xi^3 +\varphi_2 \xi^6 \right)
 \left( \varphi_9 +\varphi_6\xi^3 +\varphi_3 \xi^6 +\xi^9 \right)  \\
  \left( \varphi_{10} +\varphi_7\xi^3 +\varphi_4 \xi^6 +
          \varphi_1  \xi^9\right) \, ,
\end{multline*}
which specializes, for $\x=\{1,1,1,1,1\}$, into
$$ (1-\xi^3)^{-15} (15+51\xi^3+15\xi^6) (5+45\xi^3+30\xi^6+\xi^9)
       (1+30\xi^3+45\xi^6+5\xi^9) \, .$$

\end{document}